\documentclass[11pt, a4paper]{article}
\usepackage{amsmath,amssymb}
\usepackage{float} 
\usepackage{graphicx} 
\usepackage{subfigure}
\newtheorem{theorem}{Theorem}

\begin{document}

\title{Linearly degenerate PDEs  and  quadratic line complexes}
\author{ E.V. Ferapontov   and J. Moss
}
   \date{}
\vspace{-20mm}
   \maketitle
\vspace{-7mm}
\begin{center}
Department of Mathematical Sciences \\ Loughborough
University \\
Loughborough, Leicestershire LE11 3TU, UK \\[1ex] 
e-mails: \\
\texttt{E.V.Ferapontov@lboro.ac.uk}\\
\texttt{J.Moss@lboro.ac.uk}\\
\end{center}

\medskip

\begin{abstract}
A quadratic line complex is a three-parameter family of lines in  projective space $P^3$ specified by a single quadratic relation in the Pl\"ucker coordinates. Fixing a point ${\bf p}$ in $P^3$ and taking all lines of the complex passing through ${\bf p}$ we obtain a quadratic cone with vertex at ${\bf p}$. This family  of cones supplies $P^3$ with a conformal structure, which can be represented in the form $f_{ij}({\bf p}) dp^idp^j$ in a system of affine coordinates ${\bf p}=(p^1, p^2,  p^3)$. With this conformal structure we associate a three-dimensional second order quasilinear wave equation,
$$
\sum _{i, j}f_{ij}(u_{x_1}, u_{x_2},  u_{x_3}) u_{x_ix_j}=0,
$$
whose coefficients  can be obtained from $f_{ij}({\bf p})$ by setting $p^1=u_{x_1}, \ p^2=u_{x_2},  \  p^3=u_{x_3}$. We show that any PDE arising in this way is linearly degenerate, furthermore, any linearly degenerate PDE can be obtained by this construction. This provides a classification of linearly degenerate wave equations into eleven types, labelled by Segre symbols of the associated quadratic complexes. We classify Segre types for which the structure $f_{ij}({\bf p}) dp^idp^j$ is conformally flat, as well as Segre types for which the corresponding PDE is integrable.

\medskip


\medskip

\end{abstract}

\newpage

 \section{Introduction}
 
We study  second order  quasilinear  equations of the form
 \begin{equation}
f_{11} u_{x^1x^1} + f_{22} u_{x^2x^2} + f_{33} u_{x^3x^3} + 2 f_{12} u_{x^1x^2} + 2f_{13} u_{x^1x^3} + 2f_{23} u_{x^2x^3} =0, 
 \label{u}
 \end{equation}
where $u(x^1, x^2, x^3)$ is a function of three independent variables, and the coefficients 
$f_{ij}$ depend on  the first order derivatives  $u_{x^1}, u_{x^2}, u_{x^3}$ only. Throughout the paper we assume the nondegeneracy condition $\det f_{ij}\ne 0$. PDEs of this type, which can be called quasilinear wave equations,  
arise in a wide range of applications in mechanics, general relativity, differential geometry 
and the theory of integrable systems. The class of equations (\ref{u}) is invariant under the group ${\bf SL}(4)$ of linear  transformations of the dependent and independent variables $x^i, u$, which constitute the natural equivalence group of the problem. Transformations from the equivalence group act projectively on the space $P^3$ of first order derivatives $p^i=u_{x^i}$, and preserve conformal class of the quadratic form
\begin{equation}
f_{ij}({\bf p}) dp^idp^j.
\label{conf}
\end{equation}
This correspondence between nonlinear wave equations and conformal structures in projective space was proposed and thoroughly investigated in  \cite{B}. 

In the present paper we concentrate on the particular class of equations (\ref{u}) which are associated with quadratic complexes of lines in $P^3$.  Recall that the Pl\"ucker coordinates  of a line  through the points ${\bf p}=(p^1 : p^2 : p^3 : p^4)$ and ${\bf q}=(q^1 : q^2 : q^3 : q^4)$  are defined as   $p^{ij}= p^iq^j-p^jq^i$. They satisfy the quadratic Pl\"ucker relation, $\Omega=p^{23}p^{14}+p^{31}p^{24}+p^{12}p^{34}=0$. A quadratic line complex is a three-parameter family of lines in $P^3$ specified by an additional homogeneous quadratic relation among the Pl\"ucker coordinates,
$$
Q(p^{ij})=0.
$$
Fixing a point ${\bf p}$ in $P^3$ and taking the lines of the complex which pass through ${\bf p}$ one obtains a quadratic cone with vertex at ${\bf p}$. This family of cones supplies $P^3$ with a conformal structure. Its equation can be obtained by setting $ q^i=p^i+dp^i$ and passing to a system of affine coordinates, say,   $p^4=1, \ dp^4=0$. Expressions for the Pl\"ucker coordinates take the form $p^{4i}=dp^i, \ p^{ij}=p^idp^j-p^jdp^i, \ i, j=1, 2, 3$, and the  equation of the complex takes the so-called Monge form,
$$
Q(dp^i, \ p^idp^j-p^jdp^i)=f_{ij}({\bf p}) dp^idp^j=0.
$$
This provides the required   conformal structure (\ref{conf}), and the associated equation (\ref{u}).  Recall that the singular surface of the complex is the locus  in $P^3$ where the conformal structure (\ref{conf}) degenerates, $\det f_{ij}=0$. This is known to be Kummer's  quartic with 16 double points. It can be viewed as the locus where equation (\ref{u}) changes its type.

Quadratic line complexes have been extensively investigated in the classical works by Pl\"ucker \cite{Pl}, Kummer \cite{Kummer}, Klein \cite{Klein} and many other prominent geometers of 19-20th centuries, see also \cite{N, A} and references therein for more recent developments.  Lie \cite{Lie} studied certain classes of PDEs associated with line complexes. These included first order PDEs governing surfaces which are tangential to the cones of the associated conformal structure, and second order PDEs for surfaces whose asymptotic tangents belong to a given line complex (as well as surfaces conjugate to a given complex). Large part of this theory has nowadays become textbook material  \cite{Hudson, Jessop, SR, GH}.  We point out that the correspondence between quadratic complexes and three-dimensional nonlinear wave equations described above has not been discussed in the classical literature.   Our first result gives a characterisation of PDEs (\ref{u}) associated with quadratic complexes.

\medskip

\noindent {\it The following conditions are equivalent:

\noindent (1) Equation (\ref{u})/conformal structure (\ref{conf}) is associated with a quadratic line complex.

\noindent (2) Equation (\ref{u}) is linearly degenerate.

\noindent (3) Conformal structure (\ref{conf}) satisfies the condition
\begin{equation}
\partial_{(k}f_{ij)}=\varphi_{(k}f_{ij)},
\label{sym}
\end{equation}
here $\partial_k=\partial_{p^k}$,  $\varphi_k$ is a covector, and  brackets denote  complete symmetrization in  $i, j, k\in \{1, 2, 3\}$. 
Note that this equivalence holds in any dimesnion $\geq 3$.

}

\medskip
The equivalence of (1) and (3) is a well-known result, see e.g. \cite{Akivis, Safaryan}. The equivalence of (2) and (3) is the statement of Theorem 1, Sect. 2.3. The concept of linear degeneracy is defined in Sect. 2. Linearly degenerate PDEs are known to be quite exceptional from the point of view of solvability of the Cauchy problem: the gradient catastrophe, which is typical for genuinely nonlinear systems, does not occur,  so that one has global existence of classical solutions with sufficiently small initial data.

Based on the projective classification of quadratic line complexes by their Segre types \cite{Jessop}, we obtain a complete list of eleven normal forms of linearly degenerate PDEs of the form (\ref{u}) (Theorem 2 of Sect. 3). Recall that the Segre type of a quadratic complex is uniquely determined by the Jordan normal form of the $6\times 6$ matrix $Q\Omega ^{-1}$ where $Q$ and $\Omega$ are the symmetric matrices of the complex and the Pl\"ucker quadric, see  Sect. 3 for more details. The most generic linearly degenerate PDE corresponds to the Segre symbol $[111111]$:
$$
(a_1+a_2u_{x^3}^2+a_3u_{x^2}^2) u_{x^1x^1}+(a_2+a_1u_{x^3}^2+a_3u_{x^1}^2) u_{x^2x^2}+(a_3+a_1u_{x^2}^2+a_2u_{x^1}^2) u_{x^3x^3}+ 
$$
$$
2(\alpha u_{x^3}-a_3u_{x^1}u_{x^2}) u_{x^1x^2}+2(\beta u_{x^2}-a_2u_{x^1}u_{x^3}) u_{x^1x^3}+2(\gamma u_{x^1}-a_1u_{x^2}u_{x^3}) u_{x^2x^3}=0,
$$
here $a_i, \alpha, \beta, \gamma$  are constants such that $\alpha+\beta+\gamma=0$. The particular choice $\alpha=\beta=\gamma=0, \ a_1=a_2=a_3=1$, leads to the equation for minimal hypersurfaces in the Euclidean space $E^4$, 
$$
(1+u_{x^3}^2+u_{x^2}^2) u_{x^1x^1}+(1+u_{x^3}^2+u_{x^1}^2) u_{x^2x^2}+(1+u_{x^2}^2+u_{x^1}^2) u_{x^3x^3}+ 
$$
$$
-2u_{x^1}u_{x^2} u_{x^1x^2}-2u_{x^1}u_{x^3} u_{x^1x^3}-2u_{x^2}u_{x^3} u_{x^2x^3}=0,
$$
while the choice $ a_1=a_2=a_3=0$ results in the nonlinear wave equation, 
$$
\alpha u_{x^3} u_{x^1x^2}+\beta u_{x^2} u_{x^1x^3}+\gamma u_{x^1} u_{x^2x^3}=0,
$$
$\alpha+\beta+\gamma=0$, which  appeared  in the context of Veronese webs in 3D  \cite{Zakharevich}. 

Theorem 3 of Sect. 3 gives a complete list of complexes with the flat conformal structure (\ref{conf}) (conformal structures with vanishing Cotton tensor which is responsible for conformal flatness in three dimensions). 
Although the subject is fairly classical, we were not able to find a reference to this result:

\medskip

\noindent {\bf Theorem 3.}  {\it A quadratic  complex defines a flat conformal structure if and only if its Segre symbol is one of the following:
$$
[111(111)]^{*}, ~ [(111)(111)], ~ [(11)(11)(11)],  
$$
$$
[(11)(112)], ~ [(11)(22)], ~ [(114)], ~ [(123)], ~ [(222)], ~ [(24)], ~ [(33)].
$$
Here the asterisk denotes a particular subcase of  $[111(111)]$ where the matrix  $Q\Omega^{-1}$ has eigenvalues $(1, \epsilon, \epsilon^2, 0, 0, 0)$,   $\epsilon^3=1$.}

\medskip

Finally, in Theorem 4 of Sect. 3 we obtain a complete list of normal forms of linearly degenerate {\it integrable} equations of the form (\ref{u}). In general,  the integrability aspects of  quasilinear wave equations (\ref{u}) (not necessarily linearly degenerate) were investigated in \cite{B}, based on the method of hydrodynamic reductions \cite{Fer4}. It was shown that the moduli space of integrable  equations is $20$-dimensional. It was demonstrated in \cite{Odesskii} that coefficients of  `generic' integrable equations  can be parametrised by generalized hypergeometric functions. For linearly degenerate PDEs, the integrability  is equivalent to the existence of a linear Lax pair of the form 
$$
\psi_{x^2}=f(u_{x^1},  u_{x^2}, u_{x^3}, \lambda)  \psi_{x^1}, ~~~ 
\psi_{x^3}=g(u_{x^1},  u_{x^2}, u_{x^3},  \lambda) \psi_{x^1}, 
$$
where $\lambda$ is an auxiliary spectral parameter, so that  (\ref{u}) follows from the compatibility condition 
$ \psi_{x^2x^3}=\psi_{x^3x^2}$. It was pointed out in \cite{B} that the flatness of the conformal structure (\ref{conf}) is a necessary condition for integrability. Thus, all integrable cases are contained within the list of Theorem 3. 

\medskip

\noindent {\bf Theorem 4.} {\it A quadratic complex corresponds to an integrable PDE if and only if its Segre symbol is one of the following:
$$
 [(11)(11)(11)],  ~ [(11)(112)], ~ [(11)(22)],  ~ [(123)], ~ [(222)],  ~ [(33)].
$$
Modulo equivalence transformations (which are allowed to be complex-valued)  this leads to a complete list of normal forms of linearly degenerate integrable PDEs:

\noindent {\bf Segre symbol  $[(11)(11)(11)]$} 
$$
\alpha u_3u_{12}+\beta u_2 u_{13}+\gamma u_1u_{23}=0, ~~~ \alpha+\beta+\gamma=0,
$$
\noindent {\bf  Segre symbol $[(11)(112)]$} 
$$
u_{11}+u_1u_{23}-u_2u_{13}=0,
$$
\noindent {\bf  Segre symbol $[(11)(22)]$}
$$
u_{12}+u_2u_{13}-u_1u_{23}=0,
$$
\noindent {\bf  Segre symbol $[(123)]$}
$$
u_{22}+u_{13}+u_2 u_{33}-u_3 u_{23}=0,
$$
\noindent {\bf  Segre symbol $[(222)]$}
$$
u_{11}+u_{22}+ u_{33}=0,
$$
\noindent {\bf Segre symbol $[(33)]$} 
$$
u_{13}+u_1u_{22}-u_2u_{12}=0.
$$}
The  canonical  forms of Theorem 4  are not new: in different contexts, they have appeared before in \cite{Zakharevich, Pavlov, Shabat, Adler, Dun, Odesskii,  Manakov3, Ovsienko, Morozov}.  In particular, the same  normal forms  appeared in \cite{Odesskii} in the alternative approach to linear degeneracy based on the requirement of  `non-singular' structure of generalised Gibbons-Tsarev systems which govern hydrodynamic reductions of PDEs in question.

Sect. 4 contains remarks about  the Cauchy problem for linearly degenerate PDEs. We point out that for some linearly degenerate PDEs  (\ref{u}), the coefficients $f_{ij}$ can be represented in the form
$f_{ij}=\eta_{ij}+\varphi_{ij}$ where $\eta$ is a constant-coefficient matrix with diagonal entries $1, -1, -1$, while $\varphi_{ij}$ vanish at the `origin' $u_{x^1}=u_{x^2}=u_{x^3}=0$.  PDEs of this type can be viewed as nonlinear perturbations of the linear wave equation. Under the  so-called `null conditions' of Klainerman and Alinhac, the papers  \cite{Klainerman, Chris, John, Al1} establish  global existence of smooth solutions with  small initial data for multi-dimensional nonlinear wave equations. It remains to point out that both null conditions are automatically satisfied for linearly degenerate equations: they follow from the tensorial condition (\ref{sym}) satisfied in the vicinity of the origin. Our numerical simulations clearly demonstrate  that solutions with small initial data  do not break down, and behave essentially like solutions to the linear wave equation.

\section{Linearly degenerate PDEs}

In this section we discuss the concept of linear degeneracy for multidimensional second order PDEs.  After recalling the definition of linear degeneracy for first order quasilinear systems (Sect. 2.1), we extend it to second order quasilinear PDEs in 2D  (Sect. 2.2). In higher dimensions, the property of linear degeneracy is defined by the requirement of linear degeneracy of all travelling wave reductions of a  given PDE to two dimensions. This leads to the constraint (\ref{sym}) 
which is known to characterise conformal structures coming from quadratic line complexes (Sect. 2.3).

\subsection{Linearly degenerate first order quasilinear systems}

Let us consider a  quasilinear system
$$
{\bf v}_{t}+A({\bf v}){\bf v}_{x}=0,
$$
where ${\bf v}=(v^1, ..., v^n)$ is the vector of dependent variables,  $A$ is an $n\times n$ matrix, and  $t, x$ are  independent variables. Recall that a matrix $A$ is said to be linearly degenerate if its eigenvalues, assumed real and distinct, are constant in the direction of the corresponding  eigenvectors. Explicitly, $L_{r^i}\lambda^i=0$, no summation, where $L_{r^i}$ is  Lie derivative of the eigenvalue $\lambda^i$ in the direction of the corresponding eigenvector $r^i$. Linearly degenerate systems are quite exceptional from the point of view of solvability of the initial value problem, and have been thoroughly investigated in  literature, see e.g. \cite{R1, R2, Liu, Serre}. There exists a simple invariant criterion of linear degeneracy which does not appeal to eigenvalues/eigenvectors. 
Let us introduce the characteristic polynomial of $A$,
 $$
 det(\lambda E-A({\bf v}))={\lambda  }^n  +  f_1({\bf v}){\lambda
}^{n-1} +f_2({\bf v}){\lambda}^{n-2}+ \ldots + f_n({\bf v}).
$$
The condition of linear degeneracy can be represented in the form \cite{Fer},
$$
\nabla f_1~A^{n-1}+\nabla  f_2~A^{n-2}+\ldots  +\nabla  f_n=0,
$$
where $\nabla$ is the  gradient,  $\nabla f=({{\partial
f}\over {\partial v^1}},\ldots , {{\partial f}\over {\partial
v^n}})$, and $A^k$ denotes  $k$-th power of the matrix $A$. In the $2\times 2$ case this condition simplifies to
\begin{equation}
\nabla (tr A)~A=\nabla (det A).
\label{ldeg}
\end{equation}

\subsection{Linearly degenerate second order PDEs in 2D}

Here we consider second order equations of the form
 \begin{equation}
f_{11}(u_{t}, u_{x}) u_{tt} + 2 f_{12}(u_{t}, u_{x}) u_{tx}  +f_{22}(u_{t}, u_{x}) u_{xx}  =0. 
 \label{u1}
 \end{equation}
Setting $u_{t}=p^1, \ u_{x}=p^2$ we obtain an equivalent first order quasilinear representation,
\begin{equation}
p^1_{x}=p^2_{t}, ~~~ f_{11}(p^1, p^2)p^1_{t}+2f_{12}(p^1, p^2)p^1_{x}+f_{22}(p^1, p^2)p^2_{x}=0.
\label{2x2}
\end{equation}
We will call  PDE (\ref{u1}) {\it linearly degenerate} if this is the case for the corresponding quasilinear system (\ref{2x2}).  With ${\bf v}=(p^1, p^2)$,  the  constraint (\ref{ldeg})  leads to the conditions of linear degeneracy in the form
\begin{equation}
2\partial_1\left(\frac{f_{12}}{f_{11}}\right)+\partial_2\left(\ln\frac{f_{11}}{f_{22}}\right)=0,  ~~~ 2\partial_2\left(\frac{f_{12}}{f_{22}}\right)+\partial_1\left(\ln\frac{f_{22}}{f_{11}}\right)=0,
\label{ld2x2}
\end{equation}
$\partial_k=\partial_{p^k}$. In equivalent form, these conditions appeared in \cite{Men}, where they were solved implicitly leading to the following   result  (see also  \cite{Gvazava, Mokhov, Mukminov} for related work):

\medskip

\noindent {\bf Proposition.} { \it A generic linearly degenerate PDE of the form (\ref{u1}) can be represented in the form
\begin{equation}
u_{tt}-(v+w)u_{tx}+vwu_{xx}=0
\label{vw}
\end{equation}
where  the coefficients $v(u_t, u_x)$ and $w(u_t, u_x)$ are defined by the implicit relations
\begin{equation}
f(v)=vu_x-u_t, ~~~ g(w)=wu_x-u_t,
\label{fg}
\end{equation}
here $f, g$ are two arbitrary functions. Furthermore, by virtue of (\ref{vw}),  coefficients  $v(u_t, u_x)$ and $w(u_t, u_x)$ satisfy the equations
\begin{equation}
v_t=wv_x, ~~~ w_t=vw_x.
\label{vw1}
\end{equation}
Formulae (\ref{fg})  establish a B\"acklund transformation between the second order PDE (\ref{vw}) and the linearly degenerate  system (\ref{vw1}).}

\medskip

\centerline{\bf Proof:}

Setting in (\ref{ld2x2})
$
f_{11}=1, \ f_{12}=-(v+w)/2, \  f_{22}=vw,
$ 
one obtains a pair of uncoupled Hopf equations for $v$ and $w$: 
$\partial_2v+v\partial_1v=0, \ \partial_2w+w\partial_1w=0$. Their implicit solution leads to the  formula (\ref{fg}). Finally,   differentiating (\ref{fg}) by $t$ and $x$ one arrives at (\ref{vw1}). Notice that relations (\ref{fg}) can be rewritten in the form
$$
u_t=\frac{wf(v)-vg(w)}{v-w}, ~~~ u_x=\frac{f(v)-g(w)}{v-w},
$$
where the consistency condition,
$$
\left(\frac{f(v)-g(w)}{v-w}\right)_t=\left(\frac{wf(v)-vg(w)}{v-w}\right)_x,
$$
constitutes the general conservation law of the linearly degenerate system (\ref{vw1}). Thus,   (\ref{vw}) can be interpreted as the equation for the corresponding potential variable $u$, so that any two equations of the form (\ref{vw}) are equivalent to each other. This finishes the proof.

\medskip

\noindent {\bf Example. } For $f(v)=\sqrt {1-v^2}, \ g(w)=\sqrt {1-w^2}$ relations (\ref{fg}) reduce to one and the same quadratic equation for $v$ and $w$. Taking two different roots of this equation results in  the so-called Born-Infeld equation,
$$
(1+u_x^2)u_{tt}-2u_tu_xu_{tx}+(u_t^2-1)u_{xx}=0,
$$
which is the Euler-Lagrange equation for the area functional $\int  \sqrt{1+u_x^2-u_t^2}\ dxdt$ governing minimal surfaces in Minkowski space. 

\medskip

\noindent {\bf Remark. } Conditions (\ref{ld2x2}) can be represented in  tensorial form 
$$
\partial_{(k}f_{ij)}=\varphi_{(k}f_{ij)},
$$
here   $\varphi_k$ is a covector, and  brackets denote  complete symmetrization in the indices $i, j, k\in \{1, 2\}$. 
Explicitly, this gives
\begin{equation}
\begin{array}{c}
\partial_{1}f_{11}=\varphi_{1}f_{11}, ~~~ \partial_{2}f_{22}=\varphi_{2}f_{22}, \\
\ \\
\partial_{2}f_{11}+2\partial_{1}f_{12}=\varphi_{2}f_{11}+2\varphi_{1}f_{12}, ~~~ \partial_{1}f_{22}+2\partial_{2}f_{12}=\varphi_{1}f_{22}+2\varphi_{2}f_{12},
\end{array}
\label{lindeg}
\end{equation}
and the elimination of $\varphi_1, \varphi_2$ from the first two relations results in  (\ref{ld2x2}).

\subsection{Linearly degenerate second order PDEs in 3D and quadratic line complexes}

A three-dimensional PDE of the form (\ref{u}) 
is said to be {\it linearly degenerate} if all its traveling wave reductions to two dimensions are linearly degenerate in the sense of Sect. 2.2.
More precisely, setting
$ u(x^1, x^2, x^3)=u(\xi, \eta)+\zeta$ where $\xi, \eta, \zeta$ are arbitrary linear forms in the variables $x^i$, we obtain a two-dimensional equation of the form (\ref{u1}) for $u(\xi, \eta)$. The requirement of  linear degeneracy of {\it all}  such  reductions imposes strong constraints on the coefficients $f_{ij}$:

\begin{theorem} A PDE (\ref{u}) is linearly degenerate if and only if the corresponding  conformal structure $f_{ij}dp^idp^j$ satisfies the constraint (\ref{sym}), 
$$
\partial_{(k}f_{ij)}=\varphi_{(k}f_{ij)}.
$$
\end{theorem}

\medskip

\centerline{\bf Proof:}

\medskip

Let us seek traveling wave reductions in the form $u(x^1, x^2, x^3)=u(\xi, \eta)+\alpha x^1+\beta x^2+\gamma x^3$ where $\xi =x^1+\lambda x^3, \ \eta=x^2+\mu x^3$, and $ \alpha, \beta, \gamma, \lambda, \mu$ are arbitrary constants. We have
$$
u_{x^1}=u_{\xi}+\alpha, ~~ u_{x^2}=u_{\eta}+\beta,  ~~ u_{x^3}=\lambda u_{\xi}+\mu u_{\eta}+\gamma,
$$
as well as
$$
\begin{array}{c}
u_{x^1x^1}=u_{\xi \xi}, ~~u_{x^1x^2}=u_{\xi \eta}, ~~ u_{x^2x^2}=u_{\eta \eta},\\
 u_{x^1x^3}=\lambda u_{\xi \xi}+\mu u_{\xi \eta}, ~~ 
u_{x^2x^3}=\lambda u_{\xi \eta}+\mu u_{\eta \eta}, ~~ u_{x^3x^3}=\lambda^2 u_{\xi \xi}+2\lambda \mu u_{\xi \eta}+\mu^2 u_{\eta \eta}.
\end{array}
$$
The reduced equation (\ref{u}) takes the form
$$
a u_{\xi\xi}+2bu_{\xi\eta}+cu_{\eta\eta}=0,
$$
where
$$
a=f_{11}+2\lambda f_{13}+\lambda^2 f_{33}, ~~~ b=f_{12}+\lambda f_{23}+\mu f_{13}+\lambda \mu f_{33}, ~~~ 
c=f_{22}+2\mu f_{23}+\mu^2 f_{33},
$$
we point out that the coefficients $a, b, c$ are now viewed as functions of $u_{\xi}$ and $u_{\eta}$. For the reduced equation, conditions of linear degeneracy (\ref{lindeg}) take the form
$$
\partial_{u_{\xi}}a=\varphi_{1}a, ~~~ \partial_{u_{\eta}}c=\varphi_{2}c, ~~~
\partial_{u_{\eta}}a+2\partial_{u_{\xi}}b=\varphi_{2}a+2\varphi_{1}b, ~~~ \partial_{u_{\xi}}c+2\partial_{u_{\eta}}b=\varphi_{1}c+2\varphi_{2}b.
$$
Let us take the first condition, $\partial_{u_{\xi}}a=\varphi_{1}a$. The calculation of $\partial_{u_{\xi}}a$ gives
$$
\partial_{u_{\xi}}a=\partial_1f_{11}+\lambda \partial_3f_{11}+2\lambda(\partial_1f_{13}+\lambda\partial_3f_{13})+\lambda^2(\partial_1f_{33}+\lambda \partial_3f_{33}),
$$
which is   polynomial in $\lambda$ of degree three. We point out that, due to the presence of arbitrary constants $\alpha, \beta, \gamma$ in the expressions for $u_1, u_2, u_3$, the coefficients of this polynomial can be viewed as independent of  $\lambda, \mu$. Thus, $\varphi_1$ must be linear in $\lambda$, so that we can set
$\varphi_1\to \varphi_1+\lambda\varphi_3$ (keeping the same notation  $\varphi_1$ for the first term). Ultimately, the relation $\partial_{u_{\xi}}a=\varphi_{1}a$ takes the form
$$
\begin{array}{c}
\partial_1f_{11}+\lambda \partial_3f_{11}+2\lambda(\partial_1f_{13}+\lambda\partial_3f_{13})+\lambda^2(\partial_1f_{33}+\lambda \partial_3f_{33})=\\
\ \\
(\varphi_1+\lambda\varphi_3)(f_{11}+2\lambda f_{13}+\lambda^2 f_{33}).
\end{array}
$$
Equating terms at different powers of $\lambda$ we obtain four relations,
$$
\partial_{1}f_{11}=\varphi_{1}f_{11}, ~~~ \partial_{3}f_{33}=\varphi_{3}f_{33}, 
$$
$$
\partial_{3}f_{11}+2\partial_{1}f_{13}=\varphi_{3}f_{11}+2\varphi_{1}f_{13}, ~~~ \partial_{1}f_{33}+2\partial_{3}f_{13}=\varphi_{1}f_{33}+2\varphi_{3}f_{13}.
$$
Similar analysis of the three remaining  conditions of linear degeneracy of the reduced equation (where one should set $\varphi_2\to \varphi_2+\mu \varphi_3$) leads to the full set (\ref{sym}) of  conditions of linear degeneracy  in $3D$:
$$
\partial_{1}f_{11}=\varphi_{1}f_{11}, ~~ \partial_{2}f_{22}=\varphi_{2}f_{22}, ~~ \partial_{3}f_{33}=\varphi_{3}f_{33}, 
$$
$$
\partial_{2}f_{11}+2\partial_{1}f_{12}=\varphi_{2}f_{11}+2\varphi_{1}f_{12}, ~~ \partial_{1}f_{22}+2\partial_{2}f_{12}=\varphi_{1}f_{22}+2\varphi_{2}f_{12},
$$
$$
\partial_{3}f_{11}+2\partial_{1}f_{13}=\varphi_{3}f_{11}+2\varphi_{1}f_{13}, ~~ \partial_{1}f_{33}+2\partial_{3}f_{13}=\varphi_{1}f_{33}+2\varphi_{3}f_{13},
$$
$$
\partial_{2}f_{33}+2\partial_{3}f_{23}=\varphi_{2}f_{33}+2\varphi_{3}f_{23}, ~~ \partial_{3}f_{22}+2\partial_{2}f_{23}=\varphi_{3}f_{22}+2\varphi_{2}f_{23},
$$
$$
\partial_{1}f_{23}+\partial_{2}f_{13}+\partial_{3}f_{12}=\varphi_{1}f_{23}+\varphi_{2}f_{13}+\varphi_{3}f_{12}.
$$
On elimination of $\varphi$'s, these conditions give rise to seven first order differential constraints for $f_{ij}$. This finishes the proof of Theorem 1. 

\medskip

It is remarkable that the condition (\ref{sym}) is well-known in classical differential geometry, characterizing  conformal structures coming from quadratic line complexes, see e.g. \cite{Akivis,  Safaryan}. This provides  a general solution of the conditions of linear degeneracy (\ref{sym}), which can be described as follows. In  projective space $P^3$, the Pl\"ucker coordinates  of a line  through the points ${\bf p}=(p^1 : p^2 : p^3 : p^4)$ and ${\bf q}=(q^1 : q^2 : q^3 : q^4)$ are the six $2\times 2$ minors of the matrix
$$
\left(
\begin{array}{cccc}
p^1 & p^2 & p^3 & p^4\\
q^1 & q^2 & q^3 & q^4
\end{array}
\right),
$$
explicitly,  $p^{ij}= p^iq^j-p^jq^i$. These coordinates are known to satisfy the quadratic Pl\"ucker relation, $\Omega=p^{23}p^{14}+p^{31}p^{24}+p^{12}p^{34}=0$. A quadratic complex is a three-parameter family of lines in $P^3$ specified by an additional homogeneous quadratic relation among the Pl\"ucker coordinates,
$$
Q(p^{ij})=0.
$$
Fixing a point ${\bf p}$ in $P^3$ and taking the lines of the complex which pass through ${\bf p}$ one obtains a quadratic cone with vertex at ${\bf p}$. The family of these cones supplies $P^3$ with a conformal structure. Its equation can be obtained by setting $ q^i=p^i+dp^i$ and passing to a system of affine coordinates, say,   $p^4=1, \ dp^4=0$. The expressions for the Pl\"ucker coordinates take the form $p^{4i}=dp^i, \ p^{ij}=p^idp^j-p^jdp^i, \ i, j=1, 2, 3$, and the  equation of the complex takes the so-called Monge form,
$$
Q(dp^i, \ p^idp^j-p^jdp^i)=f_{ij}({\bf p}) dp^idp^j=0.
$$
This provides the required   conformal structure, which we always assume to be non-degenerate (that is, $\det f_{ij}$ is not identically zero). Since a homogeneous quadratic form $Q$ in six variables depends, modulo  a constant factor and a multiple of $\Omega$, on $19$ arbitrary parameters, we obtain the  $19$-parameter generic solution of the relations (\ref{sym}), that is, the $19$-parameter family of linearly degenerate PDEs (\ref{u}).  The singular  surface of the complex is the locus  in $P^3$ where the conformal structure degenerates, $\det f_{ij}=0$. This is known to be a quartic with 16 ordinary double points (Kummer's quartic). It can be viewed as a locus where the PDE (\ref{u}) changes its type.

\medskip

\noindent {\bf Example 1}.  The so-called tetrahedral complex, see \cite{Jessop}, Chapter VII, is defined by the equation $Q=b_1p^{41}p^{23}+b_2p^{42}p^{31}+b_3p^{43}p^{12}=0$. Its Monge form is 
$$
b_1dp^1(p^2dp^3-p^3dp^2)+b_2dp^2(p^3dp^1-p^1dp^3)+b_3dp^3(p^1dp^2-p^2dp^1)=0,
$$
or, equivalently, 
$$
(b_3-b_2)p^1dp^2dp^3+(b_1-b_3)p^2dp^1dp^3+(b_2-b_1)p^3dp^1dp^2=0.
$$
It corresponds to the `nonlinear wave  equation',
$$
(b_3-b_2)u_1u_{23}+(b_1-b_3)u_2u_{13}+(b_2-b_1)u_3u_{12}=0,
$$
which probably first  appeared in  \cite{Zakharevich} in the context of  Veronese webs in 3D. The associated Kummer surface degenerates into four planes: $p^1=0, \ p^2=0, \ p^3=0,$ plus the plane at infinity. The lines forming  tetrahedral complex are characterised by the property that their four points of intersection with the above planes have  constant cross-ratio (defined by $b_1, b_2, b_3$). Introducing the $1$-form
$$
\omega=(\lambda-b_2)(\lambda-b_3)u_1dx^1+(\lambda-b_1)(\lambda-b_3)u_2dx^2+(\lambda-b_1)(\lambda-b_2)u_3dx^3,
$$
which depends quadratically on the auxiliary `spectral parameter' $\lambda$, one can verify that the above PDE is equivalent to the condition $d\omega \wedge \omega=0$. Thus, the foliation $\omega=0$ is integrable for any value of $\lambda$, and  defines the structure known as a three-dimensional Veronese web.

\medskip

\noindent {\bf Example 2}. The so-called special complex, see \cite{Jessop}, Chapter VII, is defined by the equation $Q=(p^{12})^2+(p^{13})^2+(p^{23})^2-(p^{14})^2-(p^{24})^2-(p^{34})^2=0$. Its Monge form is
$$
(p^1dp^2-p^2dp^1)^2+(p^1dp^3-p^3dp^1)^2+(p^2dp^3-p^3dp^2)^2-(dp^1)^2-(dp^2)^2-(dp^3)^2=0,
$$
or, equivalently, 
$$
\begin{array}{c}
({p^2}^2+{p^3}^2-1){dp^1}^2+({p^1}^2+{p^3}^2-1){dp^2}^2+({p^1}^2+{p^2}^2-1){dp^3}^2\\
\ \\
-2p^1p^2dp^1dp^2-2p^1p^3dp^1dp^3-2p^2p^3dp^2dp^3=0.
\end{array}
$$
It corresponds to the equation
$$
\begin{array}{c}
(u_2^2+u_3^2-1)u_{11}+(u_1^2+u_3^2-1)u_{22}+(u_1^2+u_2^2-1)u_{33}-2u_1u_2u_{12}-2u_1u_3u_{13}-2u_2u_3u_{23}=0,
\end{array}
$$
which comes from the   Lagrangian  $\int \sqrt { u_1^2+u_2^2+u_3^2-1}\ dx$  governing minimal hypersurfaces $x_4=u(x_1, x_2, x_3)$ in Minkowski space with the metric $-dx_1^2-dx_2^2-dx_3^2+dx_4^2$.  
 The associated Kummer surface is the double sphere ${p^1}^2+{p^2}^2+{p^3}^2=1$. The external part of the sphere is the domain of hyperbolicity of our equation:  quadratic cones of the complex  are tangential to the sphere.  We point out that the equation for minimal surfaces is not integrable in dimensions higher than two.

\section{Normal forms of quadratic line complexes and  linearly degenerate PDEs}

In this section we utilise the  projective classification  of quadratic line complexes following \cite{Jessop}, see also \cite{Weiler}.  Complexes are characterised by the so-called Segre symbols which govern normal forms of  pairs of the associated quadratic forms. To be more precise, let $\Omega$ and $Q$ be the $6\times 6$ symmetric matrices of  quadratic forms specifying the Pl\"ucker quadric and the  complex. Then, for instance, the Segre symbol $[111111]$ means that the operator $Q\Omega^{-1}$ has purely diagonal Jordan normal form, the Segre symbol $[2 2 2]$ means that the operator $Q\Omega^{-1}$ has three $2\times 2$ Jordan blocks, etc. Since $Q$ is defined up to transformations of the form $Q\to \alpha Q+\beta \Omega$, we can always assume  $Q\Omega^{-1}$ to be  traceless. The associated conformal structures, and the corresponding linearly degenerate PDEs,  result from the equation of the complex upon setting 
$p^{ij}=p^idp^j-p^jdp^i, \ p^4=1, \ dp^4=0$, as explained in Sect. 2.3 (in some cases it will be more convenient to use  different affine projections, say, $p^1=1, \ dp^1=0$: this will be indicated explicitly where appropriate). Here is the summary of our results. Theorem 2 gives a complete list of normal forms of linearly degenerate  PDEs based on the classification of quadratic complexes (for simplicity, we use the notation $u_{x_i}=u_i, \ u_{x_i x_j}=u_{ij}$, etc). Theorem 3  provides a classification of complexes with the flat conformal structure $f_{ij}dp^idp^j$, and Theorem 4 characterises complexes corresponding to integrable PDEs. Theorems 2--4 will be proved simultaneously by going through the list of normal forms of quadratic complexes.

\begin{theorem} Any linearly degenerate PDE of the form (\ref{u}) can be brought by an equivalence transformation to one of the eleven canonical forms, labelled by  Segre symbols of the associated quadratic complexes.
\medskip

\noindent {\bf Case 1: Segre symbol  $[111111]$} 
$$
(a_1+a_2u_3^2+a_3u_2^2) u_{11}+(a_2+a_1u_3^2+a_3u_1^2) u_{22}+(a_3+a_1u_2^2+a_2u_1^2) u_{33}+ 
$$
$$
2(\alpha u_3-a_3u_1u_2) u_{12}+2(\beta u_2-a_2u_1u_3) u_{13}+2(\gamma u_1-a_1u_2u_3) u_{23}=0,
$$
$\alpha+\beta+\gamma=0$.
\medskip

\noindent {\bf Case 2: Segre symbol $[11112]$} 
$$
(\lambda u_2^2+\mu u_3^2+1) u_{11}+(\lambda u_1^2+\mu) u_{22}+(\mu u_1^2+\lambda) u_{33}+
$$
$$
2(\alpha u_3-\lambda u_1 u_2) u_{12}+2(\beta u_2-\mu u_1 u_3) u_{13}+2\gamma u_1 u_{23}=0,
$$
$\alpha+\beta+\gamma=0$.

\medskip

\noindent {\bf Case 3: Segre symbol $[1113]$}
$$
(\lambda u_2^2+\mu u_3^2+2 u_3) u_{11}+(\lambda u_1^2+\mu) u_{22}+(\mu u_1^2+\lambda) u_{33}+
$$
$$
2(\mu u_3-\lambda u_1 u_2 -1) u_{12}+2(\beta u_2 - \mu u_1 u_3 - u_1) u_{13}+ 2\gamma u_1 u_{23}=0,
$$
$\mu+\beta+\gamma=0$.

\medskip
\noindent {\bf Case 4: Segre symbol $[1122]$}
$$
(\lambda u_2^2+1) u_{11}+(\lambda u_1^2+4) u_{22}+\lambda u_{33}+
2(\alpha u_3-\lambda u_1u_2) u_{12}+2\beta  u_2 u_{13}+ 2\gamma u_1u_{23}=0,
$$
$\alpha+\beta+\gamma=0$.

\medskip

\noindent {\bf Case 5: Segre symbol $[114]$}
$$
\lambda u_{11}+(\lambda u_3^2+4) u_{22}+(\lambda u_2^2-2u_1) u_{33}+
2\alpha u_3u_{12}+2(u_3-\alpha u_2) u_{13}-2\lambda u_2u_3 u_{23}=0.
$$

\medskip

\noindent {\bf Case 6: Segre symbol $[123]$}
$$
\lambda u_{11} +(\lambda u_3^2+4) u_{22}+[\lambda u_2^2+2u_2] u_{33}+
2\alpha u_3 u_{12}+2(1-\lambda u_2)u_{13}+2(\gamma u_1-\lambda u_2u_3-u_3) u_{23}=0,
$$
$\alpha-\lambda+\gamma=0$.

\medskip

\noindent {\bf Case 7: Segre symbol $[222]$} 

\noindent Subcase 1: 
$$
u_{11}+u_{22}+u_{33}+2\alpha u_3 u_{12}+2\beta u_2 u_{13}+ 2\gamma u_1 u_{23}=0,
$$
\noindent Subcase 2: 
$$
(u_2^2+u_3^2)u_{11}+(u_1^2+u_3^2)u_{22}+(u_1^2+u_2^2)u_{33}+2(\alpha u_3-u_1u_2) u_{12}+2(\beta u_2 -u_1u_3) u_{13}+2(\gamma u_1-u_2u_3)  u_{23}=0,
$$
$\alpha+\beta+\gamma=0$.
\medskip

\noindent {\bf Case 8: Segre symbol $[15]$} 
$$
\lambda u_{11} + (\lambda u_3^2-2u_3) u_{22}+(\lambda u_2^2-4u_1) u_{33}+
2(\lambda u_3+1) u_{12}+2(2u_3-\lambda u_2) u_{13}+2(u_2-\lambda u_2u_3) u_{23}=0.
$$

\medskip

\noindent {\bf Case 9: Segre symbol $[24]$} 

\noindent Subcase 1: 
$$
u_{11}+u_{22}-2u_1 u_{33}+
2\lambda u_3 u_{12}+2(u_3-\lambda u_2) u_{13}=0.
$$
\noindent Subcase 2: 
$$
u_3^2u_{22}+(1+u_2^2)u_{33}+2 u_{12}+
2\lambda u_2 u_{13}-2(\lambda u_1+u_2u_3) u_{23}=0.
$$

\noindent {\bf Case 10: Segre symbol $[33]$} 
$$
\lambda u_{11}+(\lambda u_3^2-2u_3) u_{22}+(\lambda u_2^2-2u_2) u_{33}+
2(\lambda u_3+1) u_{12}+2(\lambda u_2+1) u_{13}-2(2\lambda u_1+\lambda u_2u_3-u_2-u_3) u_{23}=0.
$$

\medskip

\noindent {\bf Case 11: Segre symbol  $[6]$} 

\noindent Subcase 1:
$$
2u_3 u_{11}+u_{22}+2u_2 u_{33}-2u_1 u_{13}-2u_3 u_{23}=0.
$$
\noindent Subcase 2: 
$$
(u_3^2-2u_2) u_{11}-2u_3u_{22}+u_1^2 u_{33}+2u_1 u_{12}-2u_1u_3u_{13}+2u_2 u_{23}=0.
$$
\end{theorem}

Calculating the Cotton tensor (whose vanishing is responsible for conformal flatness in three dimensions)  we obtain a complete list of  quadratic complexes  with the  flat conformal structure. Recall that the flatness of $f_{ij}dp^idp^j$ is a necessary condition for integrability of the corresponding PDE \cite{B}. We observe that the requirement of conformal flatness imposes further constraints on the parameters appearing in  cases 1-11 of Theorem 2, which are characterised by certain coincidences among eigenvalues of the corresponding Jordan normal forms of $Q\Omega^{-1}$ (some Segre types do not possess conformally flat specialisations at all). In what follows we label conformally flat  subcases by their `refined' Segre symbols, e.g., the symbol
$[(11)(11)(11)]$ denotes the subcase of $[111111]$ with three pairs of coinciding eigenvalues,  the symbol
$[(111)(111)]$ denotes the subcase  with two triples of coinciding eigenvalues, etc, see \cite{Jessop}. 
Although the subject sounds very classical, we were not able to find a reference to the following result.

\begin{theorem} A quadratic  complex defines a flat conformal structure if and only if its Segre symbol is one of the following:
$$
[111(111)]^{*}, ~ [(111)(111)], ~ [(11)(11)(11)],  
$$
$$
[(11)(112)], ~ [(11)(22)], ~ [(114)], ~ [(123)], ~ [(222)], ~ [(24)], ~ [(33)].
$$
Here the asterisk denotes a particular subcase of  $[111(111)]$ where the matrix  $Q\Omega^{-1}$ has eigenvalues $(1, \epsilon, \epsilon^2, 0, 0, 0)$,   $\epsilon^3=1$. Modulo equivalence  transformations this gives the following list of normal forms of the associated PDEs:

\medskip

\noindent {\bf Segre symbol  $[111(111)]^{*}$} 
$$
(1-2u_2u_3) u_{11}+(1-2u_1u_3) u_{22}+2(u_1-u_2) u_{33}+ 
$$
$$
2(1+u_1u_3+u_2u_3) u_{12}+2(u_1u_2-u_3-u_2^2) u_{13}+2(u_1u_2+u_3-u_1^2) u_{23}=0,
$$

\medskip

\noindent {\bf Segre symbol  $[(111)(111)]$} 
$$
(u_2^2+u_3^2-1)u_{11}+(u_1^2+u_3^2-1)u_{22}+(u_1^2+u_2^2-1)u_{33}-2u_1u_2u_{12}-2u_1u_3u_{13}-2u_2u_3u_{23}=0,
$$

\medskip

\noindent {\bf Segre symbol  $[(11)(11)(11)]$} 
$$
\alpha u_3u_{12}+\beta u_2 u_{13}+\gamma u_1u_{23}=0, ~~~ \alpha+\beta+\gamma=0, 
$$

\medskip

\noindent {\bf  Segre symbol $[(11)(112)]$} 
$$
u_{11}+u_1u_{23}-u_2u_{13}=0,
$$

\medskip
\noindent {\bf  Segre symbol $[(11)(22)]$}
$$
u_{12}+u_2u_{13}-u_1u_{23}=0,
$$

\medskip
\noindent {\bf  Segre symbol $[(114)]$}
$$
u_{22}+u_1 u_{33}-u_3 u_{13}=0,
$$

\medskip

\noindent {\bf  Segre symbol $[(123)]$}
$$
u_{22}+u_{13}+u_2 u_{33}-u_3 u_{23}=0,
$$

\medskip

\noindent {\bf  Segre symbol $[(222)]$}
$$
u_{11}+u_{22}+u_{33}=0,
$$

\medskip

\noindent {\bf  Segre symbol $[(24)]$} 
$$
u_3^2u_{22}+(1+u_2^2)u_{33}+2 u_{12}-2u_2u_3 u_{23}=0,
$$

\noindent {\bf Segre symbol $[(33)]$} 
$$
u_{13}+u_1u_{22}-u_2u_{12}=0.
$$
\end{theorem}

Since conformal flatness  is a necessary condition for integrability, a complete list of linearly degenerate integrable PDEs can be obtained by going through the list of Theorem 3 and either 
calculating the integrability conditions as derived in \cite{B}, or verifying the existence of a Lax pair. A direct computation shows that the requirement of integrability eliminates Segre types $
[111(111)]^{*}, ~ [(111)(111)],  ~ [(114)],  ~ [(24)],$ leading to the following result:

\begin{theorem} A linearly degenerate PDE is integrable if and only if the corresponding complex has  one of the following Segre types:
$$
 [(11)(11)(11)],  ~ [(11)(112)], ~ [(11)(22)],  ~ [(123)], ~ [(222)],  ~ [(33)].
$$
Modulo equivalence transformations,  this leads to the  five canonical forms of linearly degenerate integrable PDEs (we exclude the linearisable  case with Segre symbol $[(222)]$).  For each integrable equation  we   present its Lax pair in the form $[X, Y]=0$ where $X$ and $Y$ are parameter-dependent vector fields which commute modulo the corresponding equation:
\medskip

\noindent {\bf Segre symbol  $[(11)(11)(11)]$} 
$$
\alpha u_3u_{12}+\beta u_2 u_{13}+\gamma u_1u_{23}=0,
$$
$\alpha+\beta+\gamma=0$. Setting $\alpha=a-b, \ \beta=b-c, \ \gamma=c-a$ we obtain the Lax pair:
$X=\partial_{x^3}-\frac{\lambda-b}{\lambda-c}\frac{u_3}{u_1}\partial_{x^1}, \ Y=\partial_{x^2}-\frac{\lambda-b}{\lambda-a}\frac{u_2}{u_1}\partial_{x^1}$.
\medskip

\noindent {\bf  Segre symbol $[(11)(112)]$} 
$$
u_{11}+u_1u_{23}-u_2u_{13}=0,
$$
Lax pair: $X=\partial_{x^1}-\lambda u_1\partial_{x^3}, \ Y=\partial_{x^2}+(\lambda^2u_1-\lambda u_2)\partial_{x^3}$.

\medskip
\noindent {\bf  Segre symbol $[(11)(22)]$}
$$
u_{12}+u_2u_{13}-u_1u_{23}=0,
$$
Lax pair: $X=\lambda\partial_{x^1}-u_1\partial_{x^3}, \ Y=(\lambda-1)\partial_{x^2}-u_2\partial_{x^3}$.

\medskip

\noindent {\bf  Segre symbol $[(123)]$}
$$
u_{22}+u_{13}+u_2 u_{33}-u_3 u_{23}=0,
$$
Lax pair: $X=\partial_{x^2}+(\lambda-u_3)\partial_{x^3}, \ Y=\partial_{x^1}+(\lambda^2-\lambda u_3+u_2)\partial_{x^3}$.

\medskip



\noindent {\bf Segre symbol $[(33)]$} 
$$
u_{13}+u_1u_{22}-u_2u_{12}=0,
$$
Lax pair: $X=\lambda \partial_{x^1}-u_1\partial_{x^2}, \ Y=\partial_{x^3}+(\lambda-u_2)\partial_{x^2}.$

\end{theorem}

\medskip

\noindent {\bf Remark 1.} The five canonical  forms from Theorem 4  are not new: in different contexts, they have appeared before in \cite{Zakharevich, Pavlov, Shabat, Adler, Dun, Odesskii,  Manakov3, Ovsienko, Morozov}.  The non-equivalence  of the above  PDEs can also  be seen by calculating the  Kummer surfaces of the corresponding line complexes. In all cases the Kummer surfaces degenerate into a collection of  planes: 

\noindent -- case 1: four planes in general position, one of them at infinity. 

\noindent -- case 2: two double planes, one of them at infinity. 

\noindent -- case 3: three planes, one of them double, with the double plane at infinity. 

\noindent -- case 4: one quadruple plane at infinity. 

\noindent -- case 5: two  planes, one of them triple, with the triple plane at infinity. 

\medskip

\noindent {\bf Remark 2.} Although all equations from Theorem 3 are not related via the equivalence group ${\bf SL}(4)$,  there may exist more complicated B\"acklund-type links between them. Thus,
let   $\alpha, \beta, \gamma$ and  $\tilde \alpha, \tilde \beta, \tilde \gamma$ be  two triplets of numbers  such that $\alpha+\beta+\gamma=0$ and $\tilde \alpha+\tilde \beta+\tilde \gamma$=0. Consider the system of two first order relations for the functions $u$ and $v$, 
$$
\alpha \tilde \gamma v_1u_3-\gamma \tilde \alpha v_3u_1=0, ~~~ \alpha
\tilde \beta v_2u_3- \beta \tilde \alpha v_3u_2=0.
$$
 Eliminating $v$ (that is, solving the above relations for $v_1$ and $v_2$ and imposing the compatibility condition $v_{12}=v_{21}$), we obtain the second order equation    $\alpha u_3u_{12}+\beta u_2u_{13}+\gamma u_1u_{23}=0$. 
Similarly, eliminating $u$ we obtain the analogous equation for $v$,  $\tilde \alpha v_3v_{12}+\tilde \beta v_2v_{13}+\tilde \gamma v_1v_{23}=0$. This construction first appeared in \cite{Zakharevich} in the context of  Veronese webs in 3D. It shows that any two integrable equations of the Segre type $[(11)(11)(11)]$ are related by a B\"acklund transformation. 
Similarly, the relations 
$$
(\lambda-1)v_2-u_2v_3=0, ~~~ \lambda v_1- u_1v_3=0
$$
provide a  B\"acklund transformation between the  equation for $u$, 
$
u_{12}+u_2u_{13}-u_1u_{23}=0,
$
and the equation for $v$, 
$
v_3v_{12}+(\lambda -1)v_2v_{13}-\lambda v_1v_{23}=0,
$
thus establishing the equivalence of integrable equations of the types  $[(11)(22)]$ and $[(11)(11)(11)]$.

\bigskip

\centerline{\bf Proof of Theorems 2--4:}

\medskip

We follow the classification of quadratic complexes as presented in \cite{Jessop}, p. 206-232. This constitutes  eleven canonical forms which are analysed case-by-case below. In each  case we calculate the conditions of vanishing of the Cotton tensor (responsible for conformal flatness in three dimensions), as well as the integrability conditions as derived in \cite{B}. Recall that conformal flatness is a necessary condition for integrability: this requirement already leads to a compact list of conformally flat subcases which  can be  checked for integrability by calculating the Lax pair. Our results are summarised as follows.

\medskip

\noindent {\bf Case 1 (generic): Segre symbol $[1 1 1 1 1 1]$.} The equation of the complex is
$$
\lambda_1(p^{12}+p^{34})^2-
\lambda_2(p^{12}-p^{34})^2+
\lambda_3 (p^{13}+p^{42})^2 -
\lambda_4(p^{13} - p^{42})^2+
\lambda_5(p^{14}+p^{23})^2-
\lambda_6(p^{14}-p^{23})^2=0,
$$
here $\lambda_i$ are the eigenvalues of $Q\Omega^{-1}$. Its Monge form  is
$$
[a_1+a_2(p^3)^2+a_3(p^2)^2] (dp^1)^2+[a_2+a_1(p^3)^2+a_3(p^1)^2] (dp^2)^2+[a_3+a_1(p^2)^2+a_2(p^1)^2] (dp^3)^2+
$$
$$
2[\alpha p^3-a_3p^1p^2] dp^1dp^2+2[\beta p^2-a_2p^1p^3] dp^1dp^3+2[\gamma p^1-a_1p^2p^3] dp^2dp^3=0,
$$
where $a_1=\lambda_5-\lambda_6, \ a_2=\lambda_3-\lambda_4, \ a_3=\lambda_1-\lambda_2, \ \alpha=\lambda_5+\lambda_6-\lambda_3-\lambda_4, \ \beta=\lambda_1+\lambda_2-\lambda_5-\lambda_6 , \ \gamma=\lambda_3+\lambda_4-\lambda_1-\lambda_2$, notice that $\alpha+\beta+\gamma=0$. The corresponding PDE takes the form
$$
(a_1+a_2u_3^2+a_3u_2^2) u_{11}+(a_2+a_1u_3^2+a_3u_1^2) u_{22}+(a_3+a_1u_2^2+a_2u_1^2) u_{33}+ 
$$
$$
2(\alpha u_3-a_3u_1u_2) u_{12}+2(\beta u_2-a_2u_1u_3) u_{13}+2(\gamma u_1-a_1u_2u_3) u_{23}=0,
$$
which is the 1st case of Theorem 2. The analysis of  integrability/conformal flatness  leads to the four subcases, depending on how many $a$'s equal zero.

\medskip

\noindent {\it Subcase 1}: $a_1=a_2=a_3=0$. This subcase, which corresponds to the so-called tetrahedral complex,  is integrable and conformally flat, leading to the nonlinear wave  equation \cite{Zakharevich},
$$
\alpha u_3u_{12}+\beta u_2u_{13}+\gamma u_1u_{23}=0.
$$
The Kummer surface of this complex consists of four planes  in $P^3$ in general position. The lines of the complex intersect these planes at four points with constant cross-ratio (depending on  $\alpha, \beta, \gamma$). The corresponding affinor  $Q\Omega^{-1}$ has three pairs of coinciding eigenvalues.
The notation for such complexes is $[(11)(11)(11)]$, see Example 1 of Sect. 2.3.

\medskip

\noindent {\it Subcase 2}: $a_1=a_2=0$. This subcase possesses no nondegenerate integrable specialisations. The conditions of conformal flatness imply $\alpha=-2\beta, \ a_3=\pm \beta$. For any choice of  the sign the corresponding affinor  $Q\Omega^{-1}$ has two triples of coinciding eigenvalues. Complexes of this type are denoted $[(111)(111)]$, and are known as `special': they consist of tangent lines to a nondegenerate quadric surface  in $P^3$. Particular example of this type is the PDE for minimal surfaces in Minkowski space, see Example 2 in Sect. 2.3.

\medskip

\noindent {\it Subcase 3}: $a_1=0$. The further analysis splits into two essentially different branches. The first branch corresponds to $\gamma=0, \ a_3=\pm a_2$, in this case we have both conformal flatness and integrability. The corresponding complexes are the same as in subcase 1, with Segre symbols $[(11)(11)(11)]$. The second branch corresponds to $\beta=\alpha, \
a_2=\pm \alpha, \ a_3^2+3\alpha^2=0$ or $\beta=\alpha, \
a_3=\pm \alpha, \ a_2^2+3\alpha^2=0$. All these subcases are conformally flat, but not integrable. They are  projectively equivalent to each other, with the same Segre  symbol $[111(111)]^*$ where the asterisk indicates that  the eigenvalues of the (traceless) operator $Q\Omega^{-1}$ are proportional to $(1, \epsilon, \epsilon^2, 0, 0, 0)$, here $\epsilon$ is a cubic root of unity, $\epsilon^3=1$. There exists an equivalent real normal form  of complexes of this type, the simplest one we found is
$$
(p^{24}+p^{14})^2+2(p^{12}+p^{34})(p^{23}+p^{31})=0.
$$ 
The corresponding Monge form is
$$
[1-2p^2p^3] (dp^1)^2+[1-2p^1p^3] (dp^2)^2+2(p^1-p^2) (dp^3)^2+
$$
$$
2[1+p^1p^3+p^2p^3] dp^1dp^2+2[p^1p^2-p^3-(p^2)^2] dp^1dp^3+2[p^1p^2+p^3-(p^1)^2] dp^2dp^3=0, 
$$
with the  associated PDE 
$$
(1-2u_2u_3) u_{11}+(1-2u_1u_3) u_{22}+2(u_1-u_2) u_{33}+ 
$$
$$
2(1+u_1u_3+u_2u_3) u_{12}+2(u_1u_2-u_3-u_2^2) u_{13}+2(u_1u_2+u_3-u_1^2) u_{23}=0,
$$
which is not integrable, although the corresponding conformal structure is flat. The associated Kummer surface is a double quadric, $2p^3+(p^1)^2-(p^2)^2$. This is the first case of Theorem 3.

\medskip

\noindent {\it Subcase 4}: all $a's$ are nonzero. Here we  have three essentially different  branches which, however, give no new examples. Thus, the first branch corresponds to $a_1=\epsilon_1 \gamma, \ a_2=\epsilon_2 \beta,  \ a_3=\epsilon_3 \alpha, \ \epsilon_i=\pm1$, 
in all these cases we have both conformal flatness and integrability. The corresponding complexes are the same as in subcase 1, with Segre symbols $[(11)(11)(11)]$. The second branch is $\alpha=\beta=\gamma=0, \  a_2=\epsilon_2 a_1, \ a_3=\epsilon_3 a_1, \ \epsilon_i=\pm1$.
This coincides with subcase 2, with Segre symbol $[(111)(111)]$. The third branch is $a_1=\epsilon_1\frac{\gamma^2}{\alpha-\beta}, \ a_2=\epsilon_2\frac{\beta^2}{\gamma -\alpha}, \ a_3=\epsilon_3\frac{\alpha^2}{\beta-\gamma}, \ \epsilon_i=\pm1,$ where $\alpha, \beta, \gamma \in \{1, \epsilon, \epsilon^2\}$ are three distinct cubic roots of unity. This is the same as subcase 3, with Segre symbol $[111(111)]^*$.

\medskip

\noindent {\bf Case 2: Segre symbol $[11112]$.} The equation of the complex is
$$
\lambda_1(p^{12}+p^{34})^2-
\lambda_2(p^{12}-p^{34})^2+
\lambda_3 (p^{13}+p^{42})^2 -
\lambda_4(p^{13} - p^{42})^2+
4 \lambda_5 p^{14}p^{23}+ (p^{14})^2=0.
$$
Its Monge form is
$$
[\lambda (p^2)^2+\mu(p^3)^2+1] (dp^1)^2+[\lambda (p^1)^2+\mu] (dp^2)^2+[\mu(p^1)^2+\lambda] (dp^3)^2+
$$
$$
2[\alpha p^3-\lambda p^1p^2] dp^1dp^2+2[\beta p^2- \mu p^1p^3] dp^1dp^3+ 2\gamma p^1 dp^2dp^3=0,
$$
where $\lambda=\lambda_1-\lambda_2$, $\mu=\lambda_3-\lambda_4$, $\alpha=-\lambda_3-\lambda_4+2\lambda_5$, $\beta=\lambda_1+\lambda_2-2\lambda_5$, $\gamma=-\alpha-\beta$, so that the corresponding PDE is 
$$
(\lambda u_2^2+\mu u_3^2+1) u_{11}+(\lambda u_1^2+\mu) u_{22}+(\mu u_1^2+\lambda) u_{33}+
$$
$$
2(\alpha u_3-\lambda u_1 u_2) u_{12}+2(\beta u_2-\mu u_1 u_3) u_{13}+2\gamma u_1 u_{23}=0.
$$
This is the 2nd case of Theorem 2. We verified that in this case  conditions of integrability are equivalent to conformal flatness, leading to the following subcases.

\noindent {\it Subcase 1}:  $\lambda=\mu=0, \ \alpha =0$ (the possibility $\lambda=\mu=0, \ \beta =0$ is equivalent to $\alpha=0$ via the interchange of indices 2 and 3), which simplifies to
$$
u_{11}+2\beta (u_2u_{13}-u_1 u_{23})=0.
$$
Modulo a rescaling this gives the corresponding subcases of Theorems 3-4. 

\noindent {\it Subcase 2}: $\beta=-\alpha, \ \lambda=\epsilon_1 \alpha, \ \mu=\epsilon_2 \alpha, \ \epsilon_i=\pm1$. One can show that subcase 2 is equivalent to subcase 1: all such complexes have the same Segre type
$[(11)(112)]$. 

\medskip

\noindent {\bf Case 3: Segre symbol $[1113]$.} The equation of the complex is
$$
\lambda_1(p^{12}+p^{34})^2-
\lambda_2(p^{12}-p^{34})^2-
\lambda_3(p^{13}-p^{42})^2+
\lambda_4(p^{13}+p^{42})^2+4\lambda_4p^{14}p^{23}+
2p^{14}(p^{13}+p^{42})=0.
$$
Its Monge form is
$$
[\lambda (p^2)^2+\mu(p^3)^2+2 p^3] (dp^1)^2+[\lambda (p^1)^2+\mu] (dp^2)^2+[\mu(p^1)^2+\lambda] (dp^3)^2+
$$
$$
2[\mu p^3 -\lambda p^1p^2-1] dp^1dp^2+2[\beta p^2 - \mu p^1p^3 - p^1] dp^1dp^3+ 2\gamma p^1 dp^2dp^3=0,
$$
where $\lambda=\lambda_1-\lambda_2$, $\mu=\lambda_4-\lambda_3$, $\beta=\lambda_1+\lambda_2-2\lambda_4$, $\gamma=-\mu-\beta$, so that the corresponding PDE is 
$$
(\lambda u_2^2+\mu u_3^2+2 u_3) u_{11}+(\lambda u_1^2+\mu) u_{22}+(\mu u_1^2+\lambda) u_{33}+
$$
$$
2(\mu u_3-\lambda u_1 u_2 -1) u_{12}+2(\beta u_2 - \mu u_1 u_3 - u_1) u_{13}+ 2\gamma u_1 u_{23}=0.
$$
This is the 3rd case of Theorem 2. One can show  that it possesses no non-degenerate integrable/conformally flat subcases.

\medskip

\noindent {\bf Case 4: Segre symbol $[1122]$.} The equation of the complex is
$$
\begin{array}{c}
\lambda_1(p^{12}+p^{34})^2-
\lambda_2(p^{12}-p^{34})^2+
4\lambda_3 p^{13}p^{42}+
4\lambda_4 p^{14 }p^{23}+(p^{13})^2+ 4(p^{23})^2=0.
\end{array}
$$
Setting $p^{ij}=p^idp^j-p^jdp^i$ and using the  affine projection $p^3=1, \ dp^3=0$ we obtain the  associated Monge equation,
$$
[\lambda (p^2)^2+1] (dp^1)^2+[\lambda (p^1)^2+4] (dp^2)^2+\lambda (dp^4)^2+
2[\alpha p^4-\lambda p^1p^2] dp^1dp^2+2\beta p^2 dp^1dp^4+ 2\gamma p^1dp^2dp^4,
$$
where $\lambda=\lambda_1-\lambda_2$, $\alpha=2\lambda_4-2\lambda_3$, $\beta=2\lambda_3-\lambda_1-\lambda_2, \ \gamma=-\alpha-\beta$, so that the corresponding PDE is 
$$
(\lambda u_2^2+1) u_{11}+(\lambda u_1^2+4) u_{22}+\lambda u_{44}+
2(\alpha u_4-\lambda u_1u_2) u_{12}+2\beta  u_2 u_{14}+ 2\gamma u_1u_{24}=0.
$$
Relabelling independent variables   gives the 4th case of Theorem 2. In this case conditions of conformal flatness are equivalent to the integrability, leading to $\lambda=\alpha=0$,
$$
u_{11}+4 u_{22}+2\beta  (u_2 u_{14}- u_1u_{24})=0.
$$
Modulo elementary changes of variables this gives the corresponding subcases of Theorems 3-4, with Segre symbol $[(11)(22)]$.

\medskip

\noindent {\bf Case 5: Segre symbol $[114]$.} The equation of the complex is
$$
\begin{array}{c}
\lambda_1(p^{12}+p^{34})^2-
\lambda_2(p^{12}-p^{34})^2+
4\lambda_3(p^{14}p^{23}+p^{42}p^{13})
+2p^{14}p^{42}+4(p^{13})^2=0.
\end{array}
$$
Setting $p^{ij}=p^idp^j-p^jdp^i$ and using the  affine projection $p^1=1, \ dp^1=0$ we obtain the  associated Monge equation,
$$
\lambda (dp^2)^2+[\lambda  (p^4)^2+4] (dp^3)^2+[\lambda (p^3)^2-2p^2] (dp^4)^2+
$$
$$
2\alpha p^4  dp^2dp^3+2[p^4-\alpha  p^3] dp^2dp^4-2\lambda p^3p^4  dp^3dp^4=0,
$$
where $\lambda=\lambda_1-\lambda_2$, $\alpha=2\lambda_3-\lambda_1-\lambda_2$, so that the corresponding PDE is 
$$
\lambda u_{22}+(\lambda u_4^2+4) u_{33}+(\lambda u_3^2-2u_2) u_{44}+
2\alpha u_4u_{23}+2(u_4-\alpha u_3) u_{24}-2\lambda u_3u_4 u_{34}=0.
$$
Relabelling independent variables   gives the 5th case of Theorem 2. One can show that this equation is not integrable.  
The condition of conformal flatness  gives $\lambda=\alpha=0$,
$$
4 u_{33}-2u_2 u_{44}+2u_4 u_{24}=0.
$$ 
Such complexes are denoted $[(114)]$. Modulo elementary changes of variables this gives the corresponding subcase of Theorem 3. 

\medskip

\noindent {\bf Case 6: Segre symbol $[123]$.} The equation of the complex is
$$
\begin{array}{c}
-\lambda_1(p^{12}-p^{34})^2+
4\lambda_2 p^{13}p^{42}+
4(p^{13})^2+
\lambda_3(4p^{14}p^{23}+(p^{12}+p^{34})^2)+2p^{14}(p^{12}+p^{34})=0.
\end{array}
$$
Setting $p^{ij}=p^idp^j-p^jdp^i$ and using the  affine projection $p^1=1, \ dp^1=0$ we obtain the  associated Monge equation,
$$
\lambda (dp^2)^2+[\lambda (p^4)^2+4] (dp^3)^2+[\lambda (p^3)^2+2p^3] (dp^4)^2+
$$
$$
2\alpha p^4 dp^2dp^3+2[1-\lambda p^3] dp^2dp^4+2[\gamma p^2-\lambda p^3p^4-p^4] dp^3dp^4=0,
$$
where $\lambda=\lambda_3-\lambda_1$, $\alpha=2\lambda_2-\lambda_1-\lambda_3, \ \gamma=\lambda-\alpha$, so that the corresponding PDE is 
$$
\lambda u_{22} +(\lambda u_4^2+4) u_{33}+(\lambda u_3^2+2u_3) u_{44}+
2\alpha u_4 u_{23}+2(1-\lambda u_3)u_{24}+2(\gamma u_2-\lambda u_3u_4-u_4) u_{34}=0.
$$
Relabelling independent variables   gives the 6th case of Theorem 2. In this case conditions of conformal flatness are equivalent to the integrability. One can show that both require  $\lambda=\alpha=\gamma=0$,  which gives
$$
2 u_{33}+u_{24}+u_3 u_{44}-u_4 u_{34}=0.
$$
Appropriate relabelings and rescalings give the corresponding subcases of Theorems 3-4, denoted $[(123)]$.

\medskip

\noindent {\bf Case 7: Segre symbol $[222]$.} Here we have two (projectively dual) subcases. In subcase 1 the equation of the complex is
$$
\begin{array}{c}
2\lambda_1 p^{12}p^{34}+
2\lambda_2 p^{13}p^{42}+
2\lambda_3 p^{14}p^{23}+
(p^{12})^2+(p^{13})^2+(p^{14})^2=0.
\end{array}
$$
Setting $p^{ij}=p^idp^j-p^jdp^i$ and using the  affine projection $p^1=1, \ dp^1=0$ we obtain the  associated Monge equation,
$$
(dp^2)^2+ (dp^3)^2+(dp^4)^2+
2\alpha p^4 dp^2dp^3+2\beta p^3 dp^2dp^4+2\gamma p^2 dp^3dp^4=0,
$$
where $\alpha=\lambda_2-\lambda_1, \ \beta =\lambda_1-\lambda_3, \  \gamma =\lambda_3-\lambda_2$,  so that the corresponding PDE is 
$$
u_{22}+u_{33}+u_{44}+2\alpha u_4 u_{23}+2\beta u_3 u_{24}+2\gamma u_2 u_{34}=0.
$$
 Setting $\alpha=\beta=\gamma=0$ we obtain the linear equation. The corresponding Segre symbol is $[(222)]$. One can show that the above PDE is not integrable/conformally flat  for nonzero values of constants. This is the linearisable subcase of Theorems 3-4.

In subcase 2 the equation of the complex is
$$
\begin{array}{c}
2\lambda_1 p^{12}p^{34}+
2\lambda_2 p^{13}p^{42}+
2\lambda_3 p^{14}p^{23}+
(p^{23})^2+(p^{24})^2+(p^{34})^2=0.
\end{array}
$$
Setting $p^{ij}=p^idp^j-p^jdp^i$ and using the  affine projection $p^1=1, \ dp^1=0$ we obtain the  associated Monge equation,
$$
((p^3)^2+(p^4)^2)(dp^2)^2+ ((p^2)^2+(p^4)^2)(dp^3)^2+((p^2)^2+(p^3)^2)(dp^4)^2+
$$
$$
2(\alpha p^4-p^2p^3) dp^2dp^3+2(\beta p^3-p^2p^4) dp^2dp^4+2(\gamma p^2 -p^3p^4) dp^3dp^4=0,
$$
so that the corresponding PDE is 
$$
(u_3^2+u_4^2)u_{22}+(u_2^2+u_4^2)u_{33}+(u_2^2+u_3^2)u_{44}+2(\alpha u_4-u_2u_3) u_{23}+2(\beta u_3 -u_2u_4) u_{24}+2(\gamma u_2-u_3u_4)  u_{34}=0.
$$
One can show that this subcase possesses no non-degenerate integrable/conformally flat specialisations (notice that for $\alpha=\beta=\gamma=0$ this PDE becomes degenerate). Relabelling independent variables   gives the  7th case of Theorem 2.

\medskip

\noindent {\bf Case 8: Segre symbol $[15]$.} The equation of the complex is
$$
\begin{array}{c}
-\lambda_1(p^{12}-p^{34})^2+
\lambda_2 (4p^{14}p^{23}+4p^{13}p^{42}+(p^{12}+p^{34})^2)
+4p^{14}p^{42}+2p^{13}(p^{12}+p^{34})=0.
\end{array}
$$
Setting $p^{ij}=p^idp^j-p^jdp^i$ and using the  affine projection $p^1=1, \ dp^1=0$ we obtain the  associated Monge equation,
$$
\lambda (dp^2)^2+[\lambda(p^4)^2 -2p^4](dp^3)^2+[\lambda(p^3)^2-4p^2] (dp^4)^2+
$$
$$
2[\lambda p^4+1] dp^2dp^3+2[2 p^4-\lambda p^3] dp^2dp^4+2[p^3 -\lambda p^3p^4] dp^3dp^4=0,
$$
where $\lambda=\lambda_2-\lambda_1$, so that the corresponding PDE is
$$
\lambda u_{22} + (\lambda u_4^2-2u_4) u_{33}+(\lambda u_3^2-4u_2) u_{44}+
2(\lambda u_4+1) u_{23}+2(2u_4-\lambda u_3) u_{24}+2(u_3-\lambda u_3u_4) u_{34}=0.
$$
One can show that this PDE possesses no  integrable/conformally flat specialisations. Relabelling independent variables   gives the  8th case of Theorem 2.

\medskip

\noindent {\bf Case 9: Segre symbol $[24]$.} Here we have two (projectively dual) subcases. In subcase 1 the equation of the complex is
$$
\begin{array}{c}
2\lambda_1p^{12}p^{34}+(p^{12})^2+
2\lambda_2 (p^{14}p^{23}+p^{13}p^{42})+2p^{14}p^{42}+(p^{13})^2=0.
\end{array}
$$
Setting $p^{ij}=p^idp^j-p^jdp^i$ and using the  affine projection $p^1=1, \ dp^1=0$ we obtain the  associated Monge equation,
$$
(dp^2)^2+(dp^3)^2 -2p^2(dp^4)^2+
2\lambda p^4 dp^2dp^3+2[p^4-\lambda p^3] dp^2dp^4=0,
$$
where $\lambda=\lambda_2-\lambda_1$, so that the corresponding PDE is
$$
u_{22}+u_{33}-2u_2 u_{44}+
2\lambda u_4 u_{23}+2(u_4-\lambda u_3) u_{24}=0.
$$
One can show that this  subcase  possesses no  integrable/conformally flat specialisations.
In  subcase 2 the equation of the complex is
$$
\begin{array}{c}
2\lambda_1p^{12}p^{34}+(p^{34})^2+
2\lambda_2 (p^{14}p^{23}+p^{13}p^{42})+2p^{13}p^{23}+(p^{42})^2=0.
\end{array}
$$
Setting $p^{ij}=p^idp^j-p^jdp^i$ and using the  affine projection $p^3=1, \ dp^3=0$ we obtain the  associated Monge equation,
$$
(p^4)^2(dp^2)^2+(1+(p^2)^2)(dp^4)^2 +
2 dp^1dp^2+2\lambda p^2 dp^1dp^4-2[\lambda p^1+p^2p^4] dp^2dp^4=0,
$$
where $\lambda=\lambda_2-\lambda_1$, so that the corresponding PDE is
$$
u_4^2u_{22}+(1+u_2^2)u_{44}+2 u_{12}+
2\lambda u_2 u_{14}-2(\lambda u_1+u_2u_4) u_{24}=0.
$$
One can show that this PDE is  not integrable, however, the corresponding conformal structure is  flat for $\lambda=0$.  This Segre type is known as $[(2 4)]$, giving the corresponding subcase of Theorem 3. 

Relabelling independent variables   gives the  9th case of Theorem 2.

\medskip

\noindent {\bf Case 10: Segre symbol $[33]$.} The equation of the complex is
$$
\begin{array}{c}
\lambda_1 (4p^{31}p^{24}+(p^{12}+p^{34})^2)+2p^{13}(p^{12}+p^{34})+
\lambda_2 (4p^{23}p^{14}-(p^{12}-p^{34})^2)+2p^{14}(p^{12}-p^{34})=0.
\end{array}
$$
Setting $p^{ij}=p^idp^j-p^jdp^i$ and using the  affine projection $p^1=1, \ dp^1=0$ we obtain the  associated Monge equation,
$$
\lambda (dp^2)^2+[\lambda (p^4)^2-2p^4] (dp^3)^2+[\lambda (p^3)^2-2p^3] (dp^4)^2+
$$
$$
2[\lambda p^4+1] dp^2dp^3+2[\lambda p^3+1] dp^2dp^4-2[2\lambda p^2+\lambda p^3p^4-p^3-p^4] dp^3dp^4=0,
$$
where $\lambda=\lambda_1-\lambda_2$, so that the corresponding PDE is
$$
\lambda u_{22}+(\lambda u_4^2-2u_4) u_{33}+(\lambda u_3^2-2u_3) u_{44}+
2(\lambda u_4+1) u_{23}+2(\lambda u_3+1) u_{24}-2(2\lambda u_2+\lambda u_3u_4-u_3-u_4) u_{34}=0.
$$
Relabelling independent variables   gives the 10th case of Theorem 2.  One can show  that the conditions of   integrability are equivalent to conformal flatness, leading to   $\lambda=0$,
$$
u_4 u_{33}+u_3 u_{44}
- u_{23}- u_{24}-(u_3+u_4) u_{34}=0.
$$
The corresponding complex is denoted $[(33)]$. Introducing the new independent variables $x, y, t$ such that $\partial_3=\partial_x+\partial_y, \ \partial_4=\partial_x-\partial_y,\ \partial_2=-2\partial_t$ one can  reduce the above PDE to the canonical form
$$
u_{xt}+u_xu_{yy}-u_yu_{xy}=0.
$$
This is the last case of Theorems 3-4.

\medskip

\noindent {\bf Case 11: Segre symbol $[6]$.} Here we have two (projectively dual) subcases. In subcase 1 the equation of the complex is
$$
\begin{array}{c}
2\lambda (p^{23}p^{14}+p^{31}p^{24}+p^{12}p^{34})+2p^{14}p^{34}+2p^{12}p^{42}+(p^{13})^2=0.
\end{array}
$$
Setting $p^{ij}=p^idp^j-p^jdp^i$ and using the  affine projection $p^1=1, \ dp^1=0$ we obtain the  associated Monge equation,
$$
2p_4(dp^2)^2+ (dp^3)^2+2p^3(dp^4)^2-2p^2 dp^2dp^4-2p^4dp^3dp^4=0,
$$
so that the corresponding PDE is
$$
2u_4 u_{22}+u_{33}+2u_3 u_{44}-2u_2 u_{24}-2u_4 u_{34}=0.
$$
In the second subcase the equation of the complex is
$$
\begin{array}{c}
2\lambda (p^{23}p^{14}+p^{31}p^{24}+p^{12}p^{34})+2p^{23}p^{12}+2p^{34}p^{13}+(p^{42})^2=0.
\end{array}
$$
Setting $p^{ij}=p^idp^j-p^jdp^i$ and using the  affine projection $p^1=1, \ dp^1=0$ we obtain the  associated Monge equation,
$$
((p^4)^2-2p_3)(dp^2)^2-2p^4 (dp^3)^2+(p^2)^2(dp^4)^2+2p^2 dp^2dp^3-2p^2p^4dp^2dp^4+2p^3dp^3dp^4=0,
$$
so that the corresponding PDE is
$$
(u_4^2-2u_3) u_{22}-2u_4u_{33}+u_2^2 u_{44}+2u_2 u_{23}-2u_2u_4u_{24}+2u_3 u_{34}=0.
$$
One can show  that both subcases are not  integrable/conformally flat. Relabelling independent variables   gives the last case of Theorem 2. 
This finished the proof of Theorems 2-4.

\section{Remarks on the Cauchy problem for linearly degenerate PDEs}

In $1+1$ dimensions, linearly degenerate systems are known to be quite exceptional from the point of view of solvability of the Cauchy problem: generic smooth initial data do not develop shocks in finite time \cite{R1, R2, Liu, Serre}. The conjecture of Majda
\cite{Majda}, p. 89, suggests that the same statement should be true in higher dimensions, namely, for linearly degenerate systems the shock formation never happens for smooth initial data. To the best of our knowledge this conjecture is largely open, and has only been established for particular classes of multi-dimensional linearly degenerate PDEs, see \cite{Klainerman, Chris, John, Brenier} and references therein. We  emphasize that the so-called `null condition' of Klainerman, which is instrumental for establishing  global existence results for $3+1$ dimensional  nonlinear wave equations with small initial data, is automatically satisfied for  linearly degenerate PDEs. In the more subtle case of $2+1$ dimensions, the null condition implies  long time existence, and additional conditions (e.g. the second null condition of Alinhac, which also follows from linear degeneracy) are required to guarantee global existence  \cite{Al1}.
The approach of \cite{Klainerman, Chris, John, Al1} applies to second order quasilinear PDEs which can be viewed as nonlinear deformations of the wave equation,
\begin{equation}
\square u= g_{ij}(u_k)u_{ij},
\label{wave}
\end{equation}
here $\square=\partial_1^2-\partial_2^2-... -\partial_n^2$ is the wave operator, and the coefficients  $g_{ij}$, which depend on the first order derivatives of $u$, are required to vanish at the origin $u_k=0$. Under the  null conditions imposed on $g_{ij}$ (which are automatically satisfied for linearly degenerate PDEs of the form (\ref{wave}), in fact, these conditions follow from the requirement of linear degeneracy  in the vicinity of  the origin), one has  global existence  for classical solutions with small initial data. Since some of the linearly degenerate examples from Theorem 2 can be put into the form (\ref{wave}), one can automatically guarantee global existence. 
For instance, the PDE for  minimal hypersurfaces in the Minkowski space is
\begin{equation}
u_{11}-u_{22}-u_{33}=-(u_2^2+u_3^2) u_{11}+(u_3^2-u_1^2) u_{22}+(u_2^2-u_1^2) u_{33}
+2u_1u_2 u_{12}+2u_1u_3 u_{13}-2u_2u_3 u_{23},
\label{C1}
\end{equation}
take case $[111111]$ of Theorem 2 and set $  a_1=-1, \ a_2=a_3=1, \ \alpha=\beta=\gamma=0, \ u\to iu$. It can be obtained as the Euler-Lagrange equation for the area functional,
$\int \sqrt {1+u_2^2+u_3^2-u_1^2} \ dx$. In this particular case  global existence was established in \cite{Lindblad}, in fact, this PDE fits into the general framework of \cite{Al1}. Further examples of this type include the equation
\begin{equation}
u_{11}-u_{22}-u_{33}=2\alpha u_3 u_{12}+2\beta u_2 u_{13}+ 2\gamma u_1 u_{23},
\label{C2}
\end{equation}
take case $[222]$ of Theorem 2 and set $ x_2\to ix_2, \ x_3\to ix_3$. For  PDEs of this type, solutions with small initial data essentially behave  like solutions of the linear wave equation. As an illustration we present Mathematica snapshots of numerical solutions for equations (\ref{C1})-(\ref{C2}) with  hump-like initial data at $x_1=0$: $u=0.8 e^{-x_2^2-x_3^2}, \ u_{x_1}=0$.

\begin{figure}[H]
\centering
\subfigure{\includegraphics[scale=0.4]{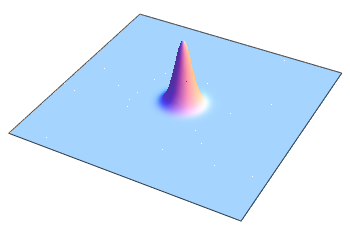}
\label{fig:subfig1}}
\subfigure{\includegraphics[scale=0.4]{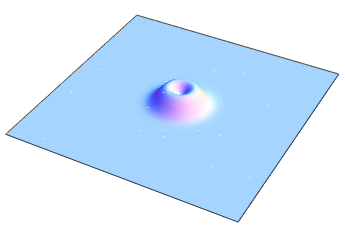}
\label{fig:subfig2}}
\subfigure{\includegraphics[scale=0.4]{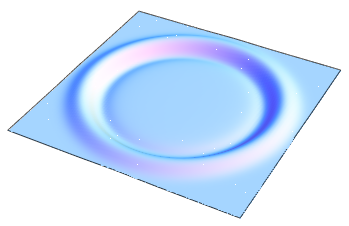}
\label{fig:subfig3}}
\label{fig:subfigureExample}
\caption{Numerical solution of equation (\ref{C1}) for  $x_1=0, 1, 8$.}
\end{figure}

\begin{figure}[H]
\centering
\subfigure{\includegraphics[scale=0.4]{21}
\label{fig:subfig1}}
\subfigure{\includegraphics[scale=0.4]{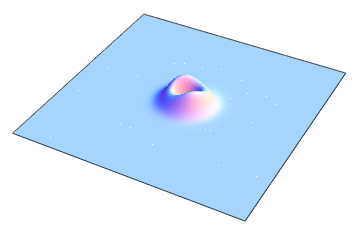}
\label{fig:subfig2}}
\subfigure{\includegraphics[scale=0.4]{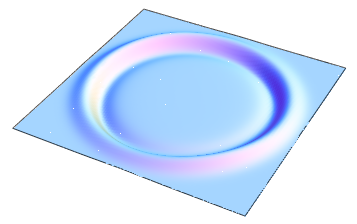}
\label{fig:subfig3}}
\label{fig:subfigureExample}
\caption{Numerical solution of equation (\ref{C2}) for $x_1=0, 1, 8$.}
\end{figure}
\noindent Although we observe some minor differences at  the early stages of evolution, for large values of $x_1$ solutions become almost indistinguishable from analogous solutions of the linear wave equation. 

We also refer to \cite{Manakov1, Manakov3} for an alternative approach to the Cauchy problem for linearly degenerate integrable PDEs based on the novel version of the inverse scattering transform.

\section*{Acknowledgements}

We thank Alexey Bolsinov,  Yann Brenier, Gavin Brown, Pavel Burovskii, Karima Khusnutdinova  and Maxim Pavlov  for clarifying discussions. The research of EVF  was partially supported by the European Research Council Advanced Grant  FroM-PDE.


\begin{thebibliography}{99}



\bibitem{Adler} V.E. Adler and  A.B. Shabat, 
Model equation of the theory of solitons, {\bf 153}, no. 1 (2007) 1373--1387.


\bibitem{Akivis} M.A. Akivis and V.V. Goldberg, Projective differential geometry of submanifolds, Elsevier Science Publishers (1993) 362pp.

\bibitem{Al1} S. Alinhac, The null condition for quasilinear wave equations in two space dimensions I, Invent. Math. {\bf 145}, no. 3 (2001) 597--618.


\bibitem{A} D. Avritzer and H. Lange, Moduli spaces of quadratic complexes and their singular surfaces, Geom. Dedicata {\bf 127} (2007) 177--197.

\bibitem{Brenier} Y. Brenier,  Hydrodynamic structure of the augmented Born-Infeld equations, Arch. Ration. Mech. Anal. {\bf 172}, no. 1 (2004) 65--91.

\bibitem{B}  P.A. Burovskii, E.V. Ferapontov and S.P. Tsarev, Second order quasilinear PDEs and conformal structures in projective space, International J. Math. { \bf 21}, no. 6 (2010) 799-841.

\bibitem{Chris} D. Christodoulou, Global solutions of nonlinear hyperbolic equations for small initial data, Comm. Pure Appl. Math. {\bf 39} (1986) 267--282.

\bibitem{Dub} B.A. Dubrovin and S.P. Novikov,
Hydrodynamics of weakly
deformed soliton lattices:
differential geometry and
Hamiltonian theory, Russian Math. Surveys,
{\bf 44}  (1989) 35--124.

\bibitem{Dun} M. Dunajski, A class of Einstein-Weyl spaces associated to an integrable system of hydrodynamic type. J. Geom. Phys. {\bf 51}, no. 1 (2004) 126--137.

\bibitem{Fer} E.V. Ferapontov, Integration of weakly nonlinear hydrodynamic
systems in Riemann invariants, Phys. Lett. A {\bf 158} (1991) 112-118. 


\bibitem{Fer4} E.V. Ferapontov and K.R. Khusnutdinova, On integrability of
(2+1)-dimensional quasilinear systems, Comm. Math. Phys.  {\bf 248} (2004)
187--206.



\bibitem{Gibb94} J. Gibbons and Y. Kodama,   A method for solving the
dispersionless KP hierarchy and its exact solutions. II. Phys. Lett.
A {\bf135} (1989) 167--170.

\bibitem{GibTsa96} J. Gibbons and S.P. Tsarev,
Reductions of the
Benney equations, Phys. Lett. A {\bf 211} (1996) 19--24.

\bibitem{GibTsa99} J. Gibbons and S.P. Tsarev, Conformal maps and
reductions of
the Benney equations, Phys. Lett. A {\bf 258} (1999)
263--271.

\bibitem{GH} P. Griffiths and J. Harris,  Principles of algebraic geometry, John Wiley and Sons, New York (1978) 813 pp. 

\bibitem{Gvazava} D.K. Gvazava, General integrals and initial-characteristic problems for second-order nonlinear equations with rectilinear characteristics,  Lithuanian Math. J. {\bf 40}, no. 4 (2000) 352--363. 

\bibitem{Hudson} R.W.H.T. Hudson, Kummer's quartic surface, Cambridge University Press, 1905. 


\bibitem{Jessop} C.M. Jessop,  A treatise on the line complex, Chelsea Publishing Co., New York (1969)  364 pp.

\bibitem{John} F. John, Existence for large times of strict solutions of nonlinear wave equations in three space dimensions for small initial data, Comm. Pure Appl. Math. {\bf 40}, no. 1 (1987) 79--109.

\bibitem{Klainerman} S. Klainerman, 
The null condition and global existence to nonlinear wave equations. Nonlinear systems of partial differential equations in applied mathematics, Part 1 (Santa Fe, N.M., 1984), 293--326, 
Lectures in Appl. Math. {\bf 23}, Amer. Math. Soc., Providence, RI, 1986. 

\bibitem{Klein} F. Klein,  Zur theorie der Liniencomplexe des ersten und zweiten Grades, Math. Ann. {\bf 2} (1870) 198-226.

\bibitem{Kummer} E. Kummer, 
 \"Uber die Fl\"achen vierten Grades mit sechzehn singul\"aren Punkten, Monatsberichte der Koniglichen Preusischen Akademie der Wissenschaften zu Berlin (1864) 246--260;
see also Collected papers, 
Volume II: Function theory, geometry and miscellaneous,  Springer-Verlag, Berlin-New York (1975) 877 pp.

\bibitem{Lie} S. Lie,  Geometrie der Ber\"uhrungstransformationen,  with editorial assistance by Georg Scheffers. Second corrected edition. Chelsea Publishing Co., New York, 1977. xi+694 pp. ISBN: 0-8284-0291-4.


\bibitem{Lindblad} H. Lindblad,  A remark on global existence for small initial data of the minimal surface equation in Minkowskian space time, Proc. Amer. Math. Soc. {\bf 132}, no.4 (2004) 1095--1102.

\bibitem{Liu} T.P. Liu, Development of singularities in the nonlinear waves for quasi-linear hyperbolic PDEs, J. Diff. Eq. {\bf 33} (1979) 92--111.	


\bibitem{Majda} A. Majda, Compressible fluid flows and systems of
conservation laws in several space variables, Appl. Math. Sci.,
Springer-Verlag, NY, {\bf 53} (1984) 159 pp.

\bibitem{Manakov1} S.V. Manakov and P.M.  Santini,  Inverse scattering problem for vector fields and the Cauchy problem for the heavenly equation, Phys. Lett. A {\bf 359}, no. 6 (2006) 613--619.




\bibitem{Manakov3} S.V. Manakov and P.M.  Santini, 
On the solutions of the second heavenly and Pavlov equations, J. Phys. A {\bf 42}, no. 40 (2009) 404013, 11 pp.


\bibitem{Shabat} L. Martinez Alonso and A.B. Shabat, Hydrodynamic reductions and
solutions of a universal hierarchy, Theoret. and Mat. Phys. {\bf 140}, no. 2 (2004) 1073--1085.

\bibitem{Men} O.F. Men'shikh,  Conservation laws and B\"acklund transformations associated with the Born-Infeld equation, Math. Notes {\bf 77}, no. 3-4 (2005) 510--522.

\bibitem{Mokhov} O.I. Mokhov and Y.  Nutku,  Bianchi transformation between the real hyperbolic Monge-Amp\`ere equation and the Born-Infeld equation. Lett. Math. Phys. {\bf 32}, no. 2 (1994) 121-123.

\bibitem{Morozov} O. I. Morozov,
Recursion Operators and Nonlocal Symmetries for Integrable rmdKP and rdDym Equations,  arXiv:1202.2308.


\bibitem{Mukminov} F.Kh. Mukminov, On the straightening of the characteristics of a second-order quasilinear equation, Theoret. and Math. Phys. {\bf 75}, no. 1 (1988) 340--345.

\bibitem{N} P.E. Newstead,  Quadratic complexes. II, Math. Proc. Cambridge Philos. Soc. {\bf 91}, no. 2 (1982) 183--206.


\bibitem{Odesskii} A. Odesskii and V. Sokolov, Integrable (2+1)-dimensional systems of hydrodynamic type, 
Theoretical and Mathematical Physics {\bf 163}, no. 2 (2010) 549--586.

\bibitem{Ovsienko} V. Ovsienko, 
Bi-Hamiltonian nature of the equation $u_{tx}=u_{xy} u_y-u_{yy} u_x$, Adv. Pure Appl. Math. {\bf 1}, no. 1 (2010) 7--17.



\bibitem{Pavlov} M.V. Pavlov,  Integrable hydrodynamic chains, J. Math. Phys. {\bf 44} (2003) 4134--4156.



\bibitem{Pl} J. Pl\"ucker, Neue Geometrie des Raumes, gegr\"undet auf die Betrachtung der geraden Linie als Raumelement, Bd. I, Teubner, Leipzig, 1868; ibid. Bd. II, 1869.






\bibitem{R1} B.L. Rozdestvenskii and A.D. Sidorenko, On the impossibility of `gradient catastrophe' for weakly nonlinear systems, Z. Vycisl. Mat. i Mat. Fiz. {\bf 7} (1967) 1176-1179.

\bibitem{R2} B.L. Rozdestvenskii and N.N. Yanenko, Systems of quasilinear equations and their applications to gas dynamics, translated from the second Russian edition by J. R. Schulenberger, Translations of Mathematical Monographs, {\bf 55} American Mathematical Society, Providence, RI (1983) 676 pp. 

\bibitem{Serre} D. Serre, Systems of conservation laws. 1. 
Hyperbolicity, entropies, shock waves, Cambridge University Press, 
(1999) 263 pp;  Systems of conservation laws. 2. 
Geometric structures, oscillations, and initial-boundary value problems, 
Cambridge University Press (2000) 269 pp. 






\bibitem{Safaryan} 
L.P. Safaryan,  Certain classes of manifolds of cones of order two in $P\sb{n}$. (Russian) Akad. Nauk Armjan. SSR Dokl. {\bf 50} (1970) 83--90. 

\bibitem{SR} J.G. Semple and L. Roth, Introduction to algebraic geometry, Oxford, at the Clarendon Press, 1949, 446 pp.



\bibitem{Tsarev} S.P. Tsarev, Geometry of Hamiltonian systems of
hydrodynamic
type. Generalized hodograph method, Izvestija AN USSR
Math. {\bf 54}  (1990)
1048--1068.

\bibitem{Weiler} A. Weiler, \"Uber die verschiedenen Gattungen der Complexe zweiten Grades, Math. Ann {\bf 7} (1873) 145--207.


\bibitem{Zakharevich} I. Zakharevich,
Nonlinear wave equation, nonlinear Riemann problem, and the twistor transform of Veronese webs, arXiv:math-ph/0006001.



\end{thebibliography}
\end{document}